\documentclass[12pt]{amsart}

\usepackage[a4paper,margin=2.9cm]{geometry}
\usepackage[T1]{fontenc}
\usepackage[utf8]{inputenc}

\usepackage{mathtools}
\usepackage{amsthm}
\usepackage{graphicx}
\usepackage{tikz-cd}
\usepackage{url}

\usepackage{newtxtext}
\usepackage[varbb]{newtxmath}
\usepackage[cal=cm, scr=boondoxo]{mathalpha}
\usepackage{dsfont}
\selectfont

\newcommand{\colim}{\mathrm{colim}}

\newcommand{\red}{\mathrm{red}}

\newcommand{\Set}{\mathrm{Set}}
\newcommand{\AffSch}{\mathrm{AffSch}}
\newcommand{\op}{\mathrm{op}}

\newcommand{\simto}{\xrightarrow{\sim}}

\DeclarePairedDelimiter{\laurentseries}{(\!(}{)\!)}

\DeclareMathOperator{\Ad}{Ad}

\DeclareMathOperator{\ArtLoc}{ArtLoc}
\DeclareMathOperator{\Aut}{Aut}

\DeclareMathOperator{\diag}{diag}

\DeclareMathOperator{\Dyn}{Dyn}

\DeclareMathOperator{\Fl}{Fl}
\DeclareMathOperator{\Gal}{Gal}
\DeclareMathOperator{\Gr}{Gr}

\DeclareMathOperator{\Hom}{Hom}
\DeclareMathOperator{\Id}{Id}
\DeclareMathOperator{\im}{im}

\DeclareMathOperator{\Res}{Res}

\DeclareMathOperator{\SHom}{\mathscr{H}\kern -2pt \textit{om}}
\DeclareMathOperator{\Spec}{Spec}

\DeclareMathOperator{\trace}{tr}

\usepackage{thmtools}
\usepackage[capitalise]{cleveref}
\providecommand{\TheoremParent}{subsection}
\providecommand{\EquationParent}{section}

\theoremstyle{definition}
\declaretheorem[parent=\TheoremParent]{Definition}

\theoremstyle{remark}
\declaretheorem[sibling=Definition]{Remark}
\theoremstyle{plain}
\declaretheorem[sibling=Definition]{Theorem}
\declaretheorem[sibling=Definition]{Lemma}
\declaretheorem[sibling=Definition]{Proposition}
\declaretheorem[sibling=Definition]{Corollary}

\crefname{Definition}{Definition}{Definitions}
\crefname{Problem}{Problem}{Problems}
\crefname{Question}{Question}{Questions}
\crefname{Remark}{Remark}{Remarks}
\crefname{Example}{Example}{Examples}
\crefname{Lemma}{Lemma}{Lemmas}
\crefname{Theorem}{Theorem}{Theorems}
\crefname{Proposition}{Proposition}{Propositions}
\crefname{Corollary}{Corollary}{Corollaries}
\crefname{Conjecture}{Conjecture}{Conjectures}
\crefname{Claim}{Claim}{Claims}

\crefname{equation}{}{}

\numberwithin{equation}{\EquationParent}

\DeclareMathOperator{\SHArtLoc}{SHArtLoc}

\title{Reducedness of twisted loop groups}
\author{Zhiyuan Ding}
\address{Institute of Mathematical Sciences, ShanghaiTech University, Shanghai, 201210, China}
\email{dingzhy@shanghaitech.edu.cn}

\dedicatory{In memory of Chang Yang}

\begin{document}
\begin{abstract}
    We give an elementary proof of the reducedness of twisted loop groups along the lines of the Kneser--Tits problem.
\end{abstract}

\maketitle

\section*{Introduction}

\subsection{The problem}
Reducedness of affine Grassmannians, affine flag varieties and loop groups has been extensively studied in the past.

Let $G_0$ be a semisimple algebraic group over a field $k$. The reducedness of its loop group $LG_0$ is equivalent to the reducedness of its affine Grassmannian $\Gr_{G_0}$. Beilinson--Drinfeld~\cite{BD95} and Laszlo--Sorger~\cite{Laszlo-Sorger} proved the reducedness of $\Gr_{G_0}$ assuming $\operatorname{char}(k)=0$. Faltings~\cite{Faltings03} extended the result to arbitrary fields $k$ assuming $G_0$ is simply connected.

The next big step was made by Pappas--Rapoport~\cite{Pappas-Rapoport08}. They made progress in two directions. In one direction, they extended the study to twisted loop groups, namely $LG$ for groups $G$ over $k\laurentseries{t}$ which are not necessarily split.
In the other direction, they found that the condition $\operatorname{char}(k)\nmid\#\pi_1(G)$ is crucial for $LG$ to be reduced, from the observation that the affine Grassmannian $\Gr_{PGL_2}$ is not reduced if $\operatorname{char}(k)=2$. Under the technical conditions that $G$ is tamely ramified and $k$ is perfect, they proved that $LG$ is geometrically reduced when $G$ is semisimple and $\operatorname{char}(k)\nmid\#\pi_1(G)$, through the equivalent statement for affine flag varieties. Fakhruddin--Haines--Lourenço--Richarz~\cite{Fakhruddin-Haines-Lourenco-Richarz25} further removed the tameness hypothesis, under an additional assumption that, if $\operatorname{char}(k)=2$, the relative root system of  $G\otimes_{k\laurentseries{t}}\widehat{k\laurentseries{t}^{\mathrm{ur}}}$ is reduced.

Recently, Lourenço~\cite{Lourenco23arXiv} used techniques in condensed mathematics to remove the tameness hypothesis when the residue field $k$ is finite. For a connected reductive algebraic group $G$ over $k\laurentseries{t}$ with $k$ being finite, he proved that its loop group $LG$ is reduced if and only if $G$ is semisimple and $\operatorname{char}(k)\nmid\#\pi_1(G)$. Future developments in condensed mathematics are expected to allow his method to generalize to arbitrary residue fields $k$.

\subsection{Main results}
The main result of this paper is an elementary proof of the following statement, using only standard results about algebraic groups and group schemes.

\begin{Theorem}[See \cref{Theorem: LG is reduced}]\label{theorem:reduced}
	Let $k$ be an algebraically closed field. Let $G$ be a connected, semisimple, simply connected, absolutely almost simple algebraic group over $k\laurentseries{t}$. Then its loop group $LG$ is reduced.
\end{Theorem}

In \cref{subsection:reduction}, we shall derive the following theorem from the above one.

\begin{Theorem}\label{theorem:geometrically reduced}
	Let $k$ be any field. Let $G$ be a connected semisimple algebraic group over $k\laurentseries{t}$, and assume the order of $\pi_1(G)$ is prime to the characteristic of $k$. Then its loop group $LG$ is geometrically reduced.
\end{Theorem}

The above theorem covers all cases in which $LG$ is expected to be geometrically reduced. No additional assumptions on the group $G$ or the residue field $k$ are required.

Results in the opposite direction have been obtained in many cases. If $G$ is a connected reductive algebraic group over $k\laurentseries{t}$, it has been proved that $LG$ is not reduced when

(i) $G$ is not semisimple and $k$ is perfect by Pappas--Rapoport~\cite[Proposition~6.5]{Pappas-Rapoport08};

(ii) $G$ is semisimple and split with $\operatorname{char}(k)\mid\#\pi_1(G)$ by Haines--Lourenço--Richarz~\cite[Proposition~7.7]{Haines-Lourenco-Richarz};

(iii) $G$ is semisimple and tamely ramified,  $k$ is perfect and  $\operatorname{char}(k)\mid\#\pi_1(G)$ by Haines--Lourenço--Richarz~\cite[Proposition~7.10]{Haines-Lourenco-Richarz};

(iv) $G$ is semisimple, $k$ is finite and $\operatorname{char}(k)\mid\#\pi_1(G)$ by Lourenço~\cite[Proposition~2.8]{Lourenco23arXiv}.

The result of Lourenço is expected to extend to arbitrary fields $k$.

\subsection{Consequences}
Since the affine Grassmannian / affine flag variety is an fpqc quotient of the loop group, we have the following two corollaries of \cref{theorem:geometrically reduced}. No perfectness assumption on the field $k$ is needed in deriving these two corollaries.

\begin{Corollary}
	Let $G_0$ be a connected semisimple algebraic group over a field $k$, and assume that the order of $\pi_1(G_0)$ is prime to the characteristic of $k$. Then its affine Grassmannian $\Gr_{G_0}$ is geometrically reduced.
\end{Corollary}

\begin{Corollary}
	Let $k$ be any field and let $G$ be a connected semisimple algebraic group over $k\laurentseries{t}$. Assume that the order of $\pi_1(G)$ is prime to the characteristic of $k$. Then for every facet $\mathbf{f}$ of the Bruhat--Tits building of $G(k\laurentseries{t})$, the affine flag variety $\Fl_{G,\mathbf{f}}$ of $G$ associated with $\mathbf{f}$ is geometrically reduced.
\end{Corollary}

\subsection{Contents of the paper}
In \cref{Section: Background on loop groups}, we review the background on loop groups. 

In \cref{Section: Reduction of the problem}, we reduce the proof of our main result \cref{theorem:reduced} to \cref{proposition:Artinian Kneser-Tits for all G}, which states that the maximal torus can be generated by unipotent radicals. It is crucial to our method that it suffices to test on Artinian local rings, as proved in \cref{test of reducedness of Artinian local rings}.

\cref{section:generation-of-the-torus} is devoted to the proof of \cref{proposition:Artinian Kneser-Tits for all G}. We decompose the problem into root groups associated with connected components of the Dynkin scheme. These root groups fall into two types. Those isotypic of type $A_1$ are Weil restrictions of $SL_2$. Those of type $A_2$ will be handled by explicit computation.

\subsection{Acknowledgements}
The author would like to thank Vladimir Drinfeld for sharing his ideas on this problem. The author is grateful to Jo\~{a}o Louren\c{c}o for valuable correspondence that clarified technical details and provided helpful suggestions. The author also thanks Zhijie Dong, Lian Duan, Thomas Haines and Daniel Skodlerack for helpful discussions.

The author thanks ChatGPT (GPT-5.6 Sol Pro) for assistance in locating references for the proof of \cref{lemma:properties-of-root-groups}. The manuscript was written entirely by the author.

\section{Background on loop groups}\label{Section: Background on loop groups}
\subsection{Ind-schemes}
In this subsection, we recall the definition and some basic properties of ind-schemes over a field.
\begin{Definition}
	Let $k$ be a field. Let $\AffSch_k$ denote the category of affine $k$-schemes. A strict ind-scheme over $k$ is a functor $\AffSch_k^{\op}\to\Set$ which admits a presentation $X\cong\colim_{i\in I}X_i$ as a filtered colimit of $k$-schemes where all transition maps $X_i\to X_j(i\le j)$ are closed immersions of $k$-schemes. The category of strict ind-schemes over $k$ is the full subcategory of functors $\AffSch_k^{\op}\to\Set$ whose objects are strict ind-schemes.
\end{Definition}

All ind-schemes in this paper will be strict ind-schemes in the above sense, and we shall usually drop the word ``strict''.

\begin{Definition}
	Let $k'/k$ be a field extension. Let $X\cong\colim_{i\in I}X_i$ be an ind-scheme over $k$. We define the base change of $X$ from $k$ to $k'$ to be $X\otimes_k k':=\colim_{i\in I}(X_i\otimes_k k')$. The definition does not depend on the choice of the presentation.
\end{Definition}

\begin{Definition}
	An ind-scheme $X$ over a field $k$ is said to be reduced if there exists a presentation $X=\colim_{i\in I}X_i$ in which each $X_i$ is reduced.
\end{Definition}

\begin{Definition}
	Let $X\cong\colim_{i\in I}X_i$ be a presentation of an ind-scheme $X$ over a field $k$. We define $X_{\red}:=\colim_{i\in I}(X_i)_{\red}$. The natural closed immersions $(X_i)_{\red}\to X_i$ induce an ind-(closed immersion) $X_{\red}\to X$.
\end{Definition}

\begin{Lemma}\label{lemma:properties-of-reduced-sub-ind-scheme}
	For an ind-scheme $X$ over a field $k$, the ind-scheme $X_{\red}$ defined above is reduced and is independent of the choice of presentations. For any reduced ind-scheme $W$ over $k$, we have a canonical bijection $\Hom(W,X_{\red})\simto\Hom(W,X)$. Moreover, $X$ is reduced if and only if the natural morphism
	$X_{\red}\to X$ is an isomorphism. \qed
\end{Lemma}

\begin{Definition}
	An ind-scheme $X$ over a field $k$ is said to be geometrically reduced if $X\otimes_k\bar{k}$ is reduced.
\end{Definition}

\begin{Lemma}\label{Lemma: Colimits of reduced ind-schemes}
	If an ind-scheme over a field $k$ is a filtered colimit of reduced ind-schemes over $k$, then it is reduced. \qed
\end{Lemma}

\begin{Lemma}\label{Lemma: Product of reduced ind-schemes}
	If $X$ and $Y$ are geometrically reduced ind-schemes over a field $k$, then the fiber product $X\times_{\Spec k}Y$ is also geometrically reduced. \qed
\end{Lemma}

\subsection{Twisted loop groups}

Let $k$ be a field. Let $K=k\laurentseries{t}$.

\begin{Definition}
	Let $X$ be a scheme over $K$. We define its loop space $LX$ to be the functor 
	\begin{align*}
		\AffSch_k^{\op}&\to\Set  \\
		\Spec R & \mapsto X(R\laurentseries{t})
	\end{align*}
\end{Definition}
When $X$ is affine of finite type over $K$, $LX$ is represented by an ind-scheme over $k$.

Let $G$ be a connected affine group over $K$.
\begin{Definition}
	We define the loop group of $G$ to be its loop space $LG$. It has a natural structure of an ind-(group scheme) over $k$.
\end{Definition}

We recall some basic properties about loop spaces and loop groups.
\begin{Lemma}\label{lemma:extension of scalars of loop groups}
	Let $k'/k$ be a field extension. Then we have an isomorphism 
	\[LG\otimes_k k' \cong L(G\otimes_{k\laurentseries{t}} k'\laurentseries{t})\]
	of ind-(group scheme)s over $k'$. \qed
\end{Lemma}

\begin{Lemma}
	Let $k'/k$ be a finite field extension. Let $G'$ be a $k'$-group scheme of finite type. Then the Weil restriction $G = \Res_{k'/k}(G')$ exists as a $k$-scheme of finite type, and it has a natural $k$-group structure.
\end{Lemma}
\begin{proof}
	The statement is a special case of \cite[Proposition~A.5.1]{Conrad-Gabber-Prasad15}.
\end{proof}

\begin{Lemma}\label{lemma:Weil restriction and loop groups}
	Let $k\laurentseries{u}$ be a finite extension of $k\laurentseries{t}$. Suppose $G=\Res_{k\laurentseries{u}/k\laurentseries{t}}H$ for some linear algebraic group $H$ over $k\laurentseries{u}$. Then we have an isomorphism 
	\[LG\cong LH\]
	of ind-(group scheme)s over $k$. \qed
\end{Lemma}

\begin{Lemma}\label{lemma:loop spaces of product of schemes}
	Let $X$ and $Y$ be affine schemes over $K$. Then we have an isomorphism
	\[L(X\times_{\Spec K}Y)\cong LX\times_{\Spec k}LY\]
	of ind-schemes over $k$. \qed
\end{Lemma}

\subsection{Reduction}\label{subsection:reduction}
In this subsection, we show that \cref{theorem:reduced} implies \cref{theorem:geometrically reduced}, following Pappas--Rapoport~\cite[6.a]{Pappas-Rapoport08}.

First, \cref{lemma:extension of scalars of loop groups} shows that the geometric reducedness of $LG$ is equivalent to the reducedness of $L(G\otimes_{k\laurentseries{t}} \bar{k}\laurentseries{t})$. Thus we may assume the field $k$ in \cref{theorem:geometrically reduced} is algebraically closed.

Next, under the assumptions that $k$ is algebraically closed, that $G$ is connected and semisimple, and that $\operatorname{char}(k)\nmid \#\pi_1(G)$, it is shown in \cite[6.a]{Pappas-Rapoport08} that one can replace $G$ by its simply connected cover.

Finally, assuming $G$ is connected, semisimple and simply connected, it is a standard fact~\cite[3.1.2]{Tits66} that $G$ is a product of Weil restrictions
\[G\cong \prod_j \Res_{K_j/K}H_j\]
where $K=k\laurentseries{t}$, each $K_j$ is a finite extension of $K$, and each $H_j$ is a semisimple, simply connected, absolutely almost simple algebraic group over $K_j$. Since $k$ is algebraically closed, we have $K_j\cong k\laurentseries{u_j}$. Then \cref{lemma:Weil restriction and loop groups,lemma:loop spaces of product of schemes} imply that 
\[LG\cong\prod_j LH_j\]
Using \cref{Lemma: Product of reduced ind-schemes}, the reduction process is completed.

\section{Reduction of the problem}\label{Section: Reduction of the problem}
In this section, we reduce the proof of our main result \cref{theorem:reduced} to \cref{proposition:Artinian Kneser-Tits for all G}. The arguments are sketched in \cite[Lemma~4.26]{Fakhruddin-Haines-Lourenco-Richarz25}. We provide more details here and formulate the setup without extra hypotheses on the group $G$.

In the rest of this paper, we shall use the following notation.

Let $k$ be an algebraically closed field. Let $K=k\laurentseries{t}$.

Let $G$ be a connected, semisimple, simply connected, absolutely almost simple algebraic group over $K$.

\subsection{Notation and conventions about Borel pairs}\label{Section: Notation and conventions about Borel pairs}
The group $G$ over $K$ is quasi-split, by a theorem of Steinberg~\cite[Theorem~1.9]{Steinberg65}, extended by Borel and Springer to the case of a non-perfect field~\cite[III, 2.3, Remark~1 after Theorem~$1'$]{Serre94CG}. 

We fix a Borel pair $T\subset B$ of $G$, where $T$ is a maximal torus and $B$ is a Borel subgroup containing $T$. Both $T$ and $B$ are defined over $K$.

We put $B^+:=B$ and let $B^-$ denote the Borel subgroup opposite to $B^+$ relative to $T$. Let $U^+$ (resp. $U^-$) denote the unipotent radical of $B^+$ (resp. $B^-$).

\subsection{The big cell and the density of rational points}

\begin{Lemma}[{\cite[Exposé~XXII, Proposition~5.6.9 and Proposition~5.9.3]{SGA3IIInew}}]\label{PBW map is an open immersion}
	The morphism  $U^-\times T\times U^+\to G$ induced by multiplication is an open immersion. \qed
\end{Lemma}

The image of $U^-\times T\times U^+$ is an open dense subscheme of $G$, called the \emph{big cell}, denoted by $C(w_0)$.

\begin{Lemma}\label{big cell contains a local ring}
	Let $A$ be a local $K$-algebra. For every morphism $f:\Spec A\to G$, there exists $g\in G(K)$ such that $f$ factors through $g\cdot C(w_0)$.
\end{Lemma}
\begin{proof}
	The statement follows from the Zariski density of $G(K)$ in $G$ \cite[Corollary~18.3]{Borel91}, together with the fact that the big cell $C(w_0)$ is dense open in $G$.
\end{proof}

\subsection{The loop elementary subgroup}
\begin{Definition}
	Let $H$ be a connected quasi-split reductive group over $K$. For every $k$-algebra $R$, we define the \emph{loop elementary subgroup}
	\[LE_H(R):=\langle R_u(B')(R\laurentseries{t})\mid\text{$B'$ is a Borel subgroup of $H$ defined over $K$}\rangle\] 
	of $LH(R)$, where $R_u(B')$ is the unipotent radical of $B'$.
\end{Definition}

\begin{Lemma}
	Let $H$ be a connected quasi-split reductive group over $K$. For $j=1,2$, let $T_j\subset B_j$ be a Borel pair of $H$, and let $U_j^+$ (resp. $U_j^-$) be the unipotent radical of $B_j$ (resp. the Borel subgroup of $H$ opposite to $B_j$ relative to $T_j$). Then we have $\langle U_1^+(A), U_1^-(A)\rangle=\langle U_2^+(A), U_2^-(A)\rangle$ for all $K$-algebras $A$.
\end{Lemma}
\begin{proof}
	For $x\in H(K)$, let $\Ad(x)$ denote the inner automorphism $g\mapsto xgx^{-1}$. According to \cite[Exposé~XXVI, Corollaire~5.2]{SGA3IIInew}, we can find $u^+\in U_1^+(K)$ and $u^-\in U_1^-(K)$ such that $B_2=\Ad(u^+u^-)B_1$. Then $\Ad(u^+u^-)T_1$ is a maximal torus of $B_2$. Hence we can find $v\in U_2^+(K)$ such that $T_2=\Ad(v)\Ad(u^+u^-)T_1$ by \cite[Exposé~XXII, Corollaire~5.6.13]{SGA3IIInew}. Since $U_2^+=\Ad(u^+u^-)U_1^+$, we see that there exists $w=vu^+u^-\in\langle U_1^+(K), U_1^-(K)\rangle$ such that $U_2^+=\Ad(w)U_1^+$, $B_2=\Ad(w)B_1$ and $T_2=\Ad(w)T_1$. The uniqueness of opposite Borel relative to a given maximal torus implies that $U_2^-=\Ad(w)U_1^-$. Thus we have $\langle U_1^+(A), U_1^-(A)\rangle\supset\langle U_2^+(A), U_2^-(A)\rangle$ for all $K$-algebras $A$. The reverse inclusion follows by the same argument.
\end{proof}

For our chosen Borel pair $T\subset B$ of the group $G$, we obtain the following corollary.
\begin{Corollary}\label{corollary:loop-elementary-group-generated-by-a-pair}
	We have $LE_G(R)=\langle LU^+(R),LU^-(R)\rangle$ for all $k$-algebras $R$.
\end{Corollary}

\subsection{Reducedness of $LU^\pm$}\label{Section: Reducedness of LU}
Let $Y_n:=(LU^+\times LU^-)^n$ for $n\ge 1$. Define an ind-scheme 
\[LU^\pm := \colim_n  Y_n\]
over $k$, where the transition map $Y_n\to Y_{n+1}$ sends $x$ to $(x,1,1)$.

For each $n$, let $\mu_n: Y_n\to LG$ denote the morphism induced by multiplication (preserving the order from left to right). They are compatible with transition maps, so we obtain a morphism of ind-schemes $\mu:LU^\pm\to LG$ over $k$.

\begin{Lemma}\label{lemma:image-of-positive-and-negative-loop-unipotent}
	We have $\im(LU^\pm(R)\xrightarrow{\mu(R)}LG(R))=LE_G(R)$ for all $k$-algebras $R$.
\end{Lemma}
\begin{proof}
	The statement follows from \cref{corollary:loop-elementary-group-generated-by-a-pair}.
\end{proof}

\begin{Proposition}[{\cite[Exposé~XXVI, Corollaire~2.5]{SGA3IIInew}}]\label{Proposition: Unipotent group isomorphic to affine space}
	As a $K$-scheme, $U^+$ is isomorphic to an affine space. \qed
\end{Proposition}

\begin{Corollary}\label{LU is reduced}
	$LU^\pm$ is reduced.
\end{Corollary}
\begin{proof}
	\cref{Proposition: Unipotent group isomorphic to affine space} implies that the ind-scheme $LU^+$ over $k$ is a filtered colimit of (infinite-dimensional) affine spaces over $k$. Hence $LU^+$ is reduced by \cref{Lemma: Colimits of reduced ind-schemes}. Similarly,  $LU^-$ is also reduced. Hence each $Y_n$ is reduced by \cref{Lemma: Product of reduced ind-schemes}. Then $LU^\pm$ is reduced by \cref{Lemma: Colimits of reduced ind-schemes}.
\end{proof}

\subsection{Test of reducedness}
Let $\ArtLoc_k$ denote the category of Artinian local $k$-algebras. Note that for an object $R\in\ArtLoc_k$, we do not assume the canonical map from $k$ to the residue field of $R$ to be an isomorphism. In particular, any field containing $k$ is an object of $\ArtLoc_k$. Let $\SHArtLoc_k$ denote the full subcategory of $\ArtLoc_k$ consisting of strictly Henselian Artinian local $k$-algebras.
\begin{Lemma}\label{test of reducedness of Artinian local rings}
	Let $X$ be an ind-(locally Noetherian) ind-scheme over a field $k$. If $X_{\red}(R) = X(R)$ holds for all $R\in\SHArtLoc_k$, then $X$ is reduced.
\end{Lemma}
\begin{proof}
	The statement follows from \cite[Lemma~8.6]{Haines-Lourenco-Richarz} together with \cite[Remark~8.7]{Haines-Lourenco-Richarz}.
\end{proof}

\subsection{Reduction to the torus}\label{section:reduction-to-the-torus}
\begin{Lemma}\label{Laurent series ring is a local ring}
	If $R\in\ArtLoc_k$, then $R\laurentseries{t}$ is a local ring. \qed
\end{Lemma}

The proof of the following statement will occupy \cref{section:generation-of-the-torus}.

\begin{Proposition}\label{proposition:Artinian Kneser-Tits for all G}
	For any $R\in\ArtLoc_k$, we have $LT(R)\subset LE_G(R)$.
\end{Proposition}

\begin{Theorem}\label{Theorem: LG is reduced}
	The loop group $LG$ is reduced.
\end{Theorem}
\begin{proof}
	Choose a facet $\mathbf{f}$ of the Bruhat--Tits building of $G(k\laurentseries{t})$. Let $X=\Fl_{G,\mathbf{f}}$ be the associated affine flag variety, which is represented by an ind-scheme of ind-(finite type) over $k$ by \cite[Theorem~1.4]{Pappas-Rapoport08}. As explained in the discussion after \cite[Theorem~6.1]{Pappas-Rapoport08}, $LG$ is reduced if and only if $X$ is. Pick $x\in X(R)$ for some $R\in\SHArtLoc_k$. By \cref{test of reducedness of Artinian local rings}, it suffices to show that $x\in X_{\red}(R)$.
	
	The quotient morphism $LG\to X$ admits sections locally for the \'etale topology~\cite[Theorem~1.4]{Pappas-Rapoport08}. Thus there exists $y\in LG(R)=G(R\laurentseries{t})$ lifting $x$. We know $R\laurentseries{t}$ is a local ring by \cref{Laurent series ring is a local ring}. Then \cref{big cell contains a local ring} implies that the morphism $y:\Spec R\laurentseries{t}\to G$ factors through $g\cdot C(w_0)$ for some $g\in G(K)$. Using \cref{PBW map is an open immersion}, we can write $y=gf_1hf_2$ where $f_1\in U^-(R\laurentseries{t}), f_2\in U^+(R\laurentseries{t})$ and  $h\in T(R\laurentseries{t})$. We have $h\in \mu(LU^\pm(R))$ by \cref{lemma:image-of-positive-and-negative-loop-unipotent,proposition:Artinian Kneser-Tits for all G}. Thus one can factor the morphism $x:\Spec R\to X$ as a composition
	\[\Spec R\to LU^\pm \xrightarrow{\mu}LG\xrightarrow{g}LG\to X.\]
	By \cref{LU is reduced}, $LU^\pm$ is reduced. Hence $x$ factors through $X_{\red}$ by \cref{lemma:properties-of-reduced-sub-ind-scheme}. 
\end{proof}

\section{Generation of the torus}\label{section:generation-of-the-torus}

This section is devoted to proving \cref{proposition:Artinian Kneser-Tits for all G}.

We recall the setting of \cref{Section: Reduction of the problem}. 

Let $k$ be an algebraically closed field. Let $K=k\laurentseries{t}$.

Let $G$ be a connected, semisimple, simply connected, absolutely almost simple algebraic group over $K$. Fix a Borel pair $T\subset B$ of $G$. Let $B^+,B^-,U^+,U^-$ be as in \cref{Section: Notation and conventions about Borel pairs}.

\subsection{Dynkin schemes}
\begin{Definition}[{\cite[Exposé~XXIV, 3.2]{SGA3IIInew}}]
	A \emph{Dynkin scheme} over a scheme $S$ is a finite twisted constant scheme $\Delta$ over $S$ equipped with a subscheme $L\subset\Delta\times_S\Delta$ disjoint from the diagonal, called the scheme of links, and a length morphism $\Delta\to\{1,2,3\}_S$.
\end{Definition}

The scheme of links decomposes each geometric fiber of a Dynkin scheme into a disjoint union of linked components. For a Dynkin scheme, it is important to distinguish the connectedness of its underlying scheme from the linkedness of its geometric diagrams.

To each Dynkin scheme $\Delta$ over $K$, we associate a Dynkin scheme $\Delta^{\min}$, which keeps the same underlying scheme and the same scheme of links, and has the length morphism defined as follows. Let $\ell$ be the original length morphism on $\Delta$. For each geometric point $x\in\Delta^{\min}$, we define $\ell^{\min}(x)$ to be $\ell(x)/d(x)$, where 
\[d(x)=\gcd\{\ell(y)\mid\text{$y$ lies in the same linked component as $x$}\}.\]
This definition descends to a morphism $\ell^{\min}:\Delta^{\min}\to\{1,2,3\}_K$ over $K$.

To each reductive group $H$ over $K$, there is an associated Dynkin scheme $\Dyn(H)$ over $K$, defined in~\cite[Exposé~XXIV, 3.3]{SGA3IIInew}. We put 
\[\mathfrak{D}:=\Dyn(G).\]
Put $I=\pi_0(\mathfrak{D})$. Let $\mathfrak{D}_i$ be the connected component of $\mathfrak{D}$ corresponding to $i\in I$, equipped with the induced structure of a Dynkin scheme. Thus each $\mathfrak{D}_i$ is connected as a $K$-scheme.


In our setting, $G$ is absolutely almost simple, so we have the following classification.

\begin{Lemma}\label{lemma:type-of-connected-component-of-dynkin-scheme}
	For each $i\in I$, exactly one of the following two cases holds:
	
	(i) The Dynkin scheme $\mathfrak{D}_i$ has empty scheme of links. In this case, $\mathfrak{D}_i$ is said to be totally unlinked, and $\mathfrak{D}_i^{\min}$ is isotypic of type $A_1$.
	
	(ii) The Dynkin scheme $\mathfrak{D}_i$ is of type $A_2$.
\end{Lemma}
\begin{proof}
	Let $K_s$ be a separable closure of $K$. The chosen Borel pair $T\subset B$ identifies $\mathfrak{D}(K_s)$ with an irreducible Dynkin diagram equipped with a continuous action of $\Gal(K_s/K)$ through diagram automorphisms. Then the statement follows from the classification of orbits of diagram automorphisms of Dynkin diagrams.
\end{proof}

\subsection{Root groups}\label{section:root-groups}
For $i\in I$, let $P_i^+$ be the parabolic subgroup of $G$ containing $B^+$ corresponding to the open and closed subscheme $\mathfrak{D}_i$ of $\mathfrak{D}$ under the bijection \cite[Exposé~XXVI, Lemme~3.8]{SGA3IIInew}. Let $P_i^-$ be the parabolic subgroup opposite to $P_i^+$ relative to $T$. Put $L_i=P_i^+\cap P_i^-$. It is a common Levi subgroup of $P_i^+$ and $P_i^-$. Let $H_i$ be the derived subgroup of $L_i$. Then $H_i$ is semisimple.

We put $B_i^\pm=H_i\cap B^\pm,U_i^\pm=H_i\cap U^\pm$ and $T_i=H_i\cap T$. Applying \cite[Exposé~XXII, Corollaire~5.10.2 and Proposition~6.2.8]{SGA3IIInew}, we see that $B_i^+, B_i^-$ are opposite  Borel subgroups of $H_i$ relative to the maximal torus $T_i$, and $U_i^+$ (resp. $U_i^-$) is the unipotent radical of $B_i^+$ (resp. $B_i^-$). In particular, every $H_i$ is quasi-split.

\begin{Lemma}\label{lemma:properties-of-root-groups}
	For every $i\in I$, the group $H_i$ is simply connected, and we have a canonical isomorphism $\Dyn(H_i)\cong\mathfrak{D}_i^{\min}$ of Dynkin schemes over $K$. The morphism $\prod_{i\in I} T_i \to T$ induced by multiplication is an isomorphism.
\end{Lemma}
\begin{proof}
	It suffices to prove the assertions after base change to a separable closure $K_s$ of $K$. For a subgroup $H$ of $G$, we put $\overline{H}:=H\otimes_K K_s$. Let $\mathcal{R}=(X,\Phi,\Delta,X^\vee,\Phi^\vee,\Delta^\vee)$ be the based root datum for $(\overline{G},\overline{B^+},\overline{T})$. We have a partition $\Delta=\coprod_{i\in I}\Delta_i$ of simple roots given by  connected components $I=\pi_0(\mathfrak{D})$.
	
	From the definition of $L_i=P_i^+\cap P_i^-$ we see that $\overline{L_i}$ is the Levi subgroup of $\overline{G}$ associated with the subset $\Delta_i \subset\Delta$. Then \cite[Exposé~XXII, 5.10.7 and Th\'eor\`em~1.1]{SGA3IIInew} shows that the root datum for $(\overline{L_i},\overline{B^+\cap L_i},\overline{T})$ is $\mathcal{R}_{\Delta_i}=(X,\Phi\cap\mathbb{Z}\Delta_i,\Delta_i,X^\vee,\Phi^\vee\cap\mathbb{Z}\Delta_i^\vee,\Delta_i^\vee)$. Hence by \cite[Exposé~XXII, Propositions~6.2.7 and 6.2.8]{SGA3IIInew} the root datum for $(\overline{H_i},\overline{B_i^+},\overline{T_i})$ is the associated derived root datum of $\mathcal{R}_{\Delta_i}$, denoted by $\operatorname{d\text{\'{e}}r}(\mathcal{R}_{\Delta_i})$ in \cite[Exposé~XXI]{SGA3IIInew}. Since we assumed $G$ to be simply connected, it follows from \cite[Exposé~XXI, Lemme~6.5.11]{SGA3IIInew} that each $\overline{H_i}$ is simply connected. Since the operation $\mathcal{R}_{\Delta_i}\mapsto \operatorname{d\text{\'{e}}r}(\mathcal{R}_{\Delta_i})$ preserves the associated Dynkin diagrams, and the identification of simple roots is compatible with the Galois action, we obtain a canonical isomorphism $\Dyn(H_i)\cong\mathfrak{D}_i^{\min}$ of Dynkin schemes over $K$ by descent.
	
	Now the cocharacter lattice of $\overline{H_i}$ (resp. $\overline{G}$) is generated by $\Delta_i^\vee$ (resp. $\Delta^\vee$). The partition $\Delta=\coprod_{i\in I}\Delta_i$ of simple roots induces a partition of simple coroots $\Delta^\vee=\coprod_{i\in I}\Delta_i^\vee$ via the bijection $\Delta\to\Delta^\vee$, which shows that the base change of the morphism $\prod_{i\in I}T_i\to T$ to $K_s$ is an isomorphism $(\prod_{i\in I} T_i)\otimes_K K_s \simto T\otimes_K K_s$.
\end{proof}

\subsection{The group $SL_2$}\label{Section: SL2}
We define the following subgroups of $SL_2$ 
\[T_{SL_2}:=\begin{pmatrix}
	*  & 0 \\
	0 & *
\end{pmatrix},
B_{SL_2}^+:=\begin{pmatrix}
	* & * \\
	0 & *
\end{pmatrix},
B_{SL_2}^-:=\begin{pmatrix}
	* & 0 \\
	* & *
\end{pmatrix},
U_{SL_2}^+:=\begin{pmatrix}
	1 & * \\
	0 & 1
\end{pmatrix},
U_{SL_2}^-:=\begin{pmatrix}
	1 & 0 \\
	* & 1
\end{pmatrix}.
\]

\begin{Lemma}\label{Artinian Kneser-Tits for SL_2}
	For any ring $A$, we have $T_{SL_2}(A)\subset 
	\langle U_{SL_2}^+(A),U_{SL_2}^-(A)\rangle$.
\end{Lemma}
\begin{proof}
	For $a\in A^\times$, we have the identity
	\[
	\begin{pmatrix}
		a & 0 \\
		0 & a^{-1}
	\end{pmatrix} = 
	\begin{pmatrix}
		1 & a \\
		0 & 1
	\end{pmatrix}
	\begin{pmatrix}
		1 & 0 \\
		-a^{-1} & 1
	\end{pmatrix}
	\begin{pmatrix}
		1 & -1 \\
		0 & 1
	\end{pmatrix}
	\begin{pmatrix}
		1 & 0 \\
		1 & 1
	\end{pmatrix}
	\begin{pmatrix}
		1 & -1 \\
		0 & 1
	\end{pmatrix}
	\begin{pmatrix}
		1 & 0 \\
		-a & 1
	\end{pmatrix}. \qedhere\]
\end{proof}

\begin{Remark}
	It is easy to prove that $SL_2(A)$ is generated by $U_{SL_2}^+(A)$ and $U_{SL_2}^-(A)$ for all local rings $A$, but the statement is false for general $A$. See \cite[Section~4.3B]{Hahn-O'Meara} for discussions about counterexamples.
\end{Remark}

\begin{Lemma}\label{Artinian Kneser-Tits for Weil restriction of SL_2}
	Let $L$ be a finite separable extension of $K$. We define $T'=\Res_{L/K}T_{SL_2}$, $U'^+=\Res_{L/K}U^+_{SL_2}$, $U'^-=\Res_{L/K}U^-_{SL_2}$. Then we have $LT'(R)\subset \langle LU'^+(R),LU'^-(R)\rangle$ for every $k$-algebra $R$.
\end{Lemma}
\begin{proof}
	Let $A=R\laurentseries{t}\otimes_K L$. We have $LT'(R)=T_{SL_2}(A)$, $LU'^+(R)=U_{SL_2}^+(A)$, $LU'^-(R)=U_{SL_2}^-(A)$. Now the statement follows from \cref{Artinian Kneser-Tits for SL_2}.
\end{proof}

\subsection{Isotypic of type $A_1$}\label{section:isotypic-of-type-A_1}

Let $i\in I$ and let $\mathfrak{D}_i=\Spec K_i$ where $K_i$ is a finite separable extension of $K$. Let $H_i,B_i^+,B_i^-,U_i^+,U_i^-$ be as in \cref{section:root-groups}.

\begin{Lemma}\label{lemma:classification-isotypic-A_1}
	If $\mathfrak{D}_i$ is totally unlinked, then there is an isomorphism $H_i\simto\Res_{K_i/K}SL_2$ which maps the Borel pair $T_i\subset B_i^+$ to the Borel pair $\Res_{K_i/K}T_{SL_2}\subset\Res_{K_i/K}B_{SL_2}^+$.
\end{Lemma}
\begin{proof}
	\cref{lemma:properties-of-root-groups} shows that $H_i$ is isotypic of type $A_1$. Then \cite[Exposé~XXIV, Proposition 5.9]{SGA3IIInew} gives an isomorphism $H_i\simto\Res_{K_i/K}H_i'$ where the group $H_i'$ over $K_i$ is semisimple, simply connected and almost simple of type $A_1$. Since $H_i$ is quasi-split, $H_i'$ is also quasi-split~\cite[Exposé~XXIV, Corollaire 5.14]{SGA3IIInew}. Hence $H_i'\cong SL_2$. It remains to use the fact that all Borel pairs are conjugate~\cite[Exposé~XXVI, Corollaire 5.5]{SGA3IIInew}.
\end{proof}

\begin{Lemma}\label{lemma:Kneser-Tits-for-isotypic-of-type-A_1}
	If $\mathfrak{D}_i$ is totally unlinked, then $LT_i(R)\subset LE_{H_i}(R)$ for every $k$-algebra $R$.
\end{Lemma}
\begin{proof}
	The statement follows from \cref{Artinian Kneser-Tits for Weil restriction of SL_2,lemma:classification-isotypic-A_1}.
\end{proof}

\subsection{Kneser--Tits for $SU_3$}\label{section:Kneser-Tits-for-SU3}
Let $L/K$ be a quadratic Galois extension and let $\sigma$ be the non-trivial element of $\Gal(L/K)$. For a $K$-algebra $A$, we also write $\sigma$ for the induced automorphism $\sigma\otimes_K \Id_A \in \Aut(L\otimes_K A)$.

We define the quasi-split algebraic group $SU_3$ by
\[SU_3(A)=\{g\in SL_3(L\otimes_K A) \mid J\sigma(g^T)^{-1}J=g\}\]
for every $K$-algebra $A$, where the matrix $J$ is 
\[J=\begin{pmatrix}
	0 & 0 & 1 \\
	0 & -1 & 0 \\
	1 & 0 & 0
\end{pmatrix}.\]
We have $J^2=1$.

We take $B_{SU_3}^+$ (resp. $B_{SU_3}^-$) to be the subgroup of $SU_3$ consisting of upper (resp. lower) triangular matrices. Then they are opposite Borels relative to the maximal torus 
\[T_{SU_3}=B_{SU_3}^+\cap B_{SU_3}^-=\left\{ \begin{pmatrix}
	y & 0 & 0 \\
	0 & y^{-1}\sigma(y) & 0 \\
	0 & 0 & \sigma(y)^{-1}
\end{pmatrix}\right\}.\]
Let $U_{SU_3}^+$ (resp. $U_{SU_3}^-$) be the unipotent radical of $B_{SU_3}^+$ (resp. $B_{SU_3}^-$). We have an explicit description
\[U_{SU_3}^+=\left\{ \begin{pmatrix}
	1 & u & w \\
	0 & 1 & \sigma(u) \\
	0 & 0 & 1
\end{pmatrix} \,\middle|\, w+\sigma(w)=u\sigma(u)\right\},\]
and $U_{SU_3}^-$ is the transpose of $U_{SU_3}^+$.

The Kneser--Tits conjecture holds for quasi-split semisimple simply connected groups (see \cite[Corollaire~5.10 and Théorème~6.1]{Gille09} and \cite[Proposition~6.2(v,vi)]{Borel-Tits}). There is a quick proof~\cite{CasselmanSU3} in the special case $G=SU_3$. We include the proof below for completeness.

\begin{Proposition}\label{Field Kneser-Tits for SU_3}
	For any field extension $F/K$, we have $T_{SU_3}(F)\subset \langle U_{SU_3}^+(F), U_{SU_3}^-(F)\rangle$.
\end{Proposition}
\begin{proof}
	Let $A=L\otimes_K F$. The element $\sigma\in\Gal(L/K)$ acts on $A$ through $L$. Define the trace map $\trace: A \to F$ sending $x$ to $x+\sigma(x)$.
	For $\lambda\in F^\times$ and $y\in A^\times$, define
	\[a_\lambda=\begin{pmatrix}
		\lambda & 0 & 0 \\
		0 & 1 & 0  \\
		0 & 0 & \frac{1}{\lambda}
	\end{pmatrix},
	\quad
	b_y=\begin{pmatrix}
		y & 0 & 0 \\
		0 & \frac{\sigma(y)}{y} & 0  \\
		0 & 0 & \frac{1}{\sigma(y)}
	\end{pmatrix},
	\quad
	n_y=\begin{pmatrix}
		0 & 0 & y \\
		0 & -\frac{\sigma(y)}{y} & 0  \\
		\frac{1}{\sigma(y)} & 0 & 0
	\end{pmatrix}\]	
	to be elements of $SU_3(F)$.
	Then elements of $T_{SU_3}(F)$ are of the form $b_y$ for some $y\in A^\times$.
	
	Let $E(F)$ denote the subgroup of $SU_3(F)$ generated by $U_{SU_3}^+(F)$ and $U_{SU_3}^-(F)$.
	
	If $y\in A^\times$ satisfies $\trace(y)=x\sigma(x)$ for some $x\in A$, we have
	\[
	n_y=
	\begin{pmatrix}
		1 & x & y \\
		0 & 1 & \sigma(x)  \\
		0 & 0 & 1
	\end{pmatrix}
	\begin{pmatrix}
		1 & 0 & 0 \\
		-\frac{\sigma(x)}{\sigma(y)} & 1 & 0  \\
		\frac{1}{\sigma(y)} & -\frac{x}{y} & 1
	\end{pmatrix}
	\begin{pmatrix}
		1 & \frac{\sigma(y)x}{y} & y \\
		0 & 1 & \frac{y\sigma(x)}{\sigma(y)} \\
		0 & 0 & 1
	\end{pmatrix}\]
	where all three matrices on the right-hand side are contained in $SU_3(F)$. 
	In particular, we have $n_y\in E(F)$ for all $y\in A^\times$ with $\trace(y)=0$.	By \cref{Lemma: Trace map is surjective}, there exists $z\in L$  such that $z\ne 0$ and  $\trace(z)=0$. Then for every $\lambda\in F^\times$, we have $\trace(z\otimes\lambda)=0$ and $n_{z\otimes\lambda}n_{-z}=a_\lambda$. This shows that $a_\lambda\in E(F)$ for all $\lambda\in F^\times$. For $y\in A^\times$ with $\trace(y)\ne 0$, we have $\trace(y)\in F^\times$. Now $y'=y\cdot \trace(y)$ satisfies $\trace(y')=\trace(y)^2=\trace(y)\sigma(\trace(y))$. We have $n_y=n_{y\cdot\trace(y)}a_{\trace(y)}\in E(F)$. Therefore, we have $n_y\in E(F)$ for all $y\in A^\times$. Thus $b_y=n_y n_1\in E(F)$ for all $y\in A^\times$.
\end{proof}

\begin{Lemma}\label{Lemma: Trace map is surjective}
	For any finite separable extension $K'/K$, the trace map $K'\to K$ is surjective and $K$-linear. \qed
\end{Lemma}

\subsection{Infinitesimal lifting of Kneser--Tits for $SU_3$}
In this subsection, we prove $LT_{SU_3}(R)\subset LE_{SU_3}(R)$  for every $R\in\ArtLoc_k$. Related torus-generation calculations and infinitesimal-lifting arguments appear in the proofs of \cite[Proposition~9.3 and Corollary~9.5]{Pappas-Rapoport08}, under the tameness hypothesis.

\begin{Definition}
	A surjective ring homomorphism $R\to R_0$ is said to be a \emph{square-zero extension} by an ideal $\mathfrak{a}\subset R$ if $\mathfrak{a}=\ker(R\to R_0)$ and $\mathfrak{a}^2=0$.
\end{Definition}

\begin{Lemma}\label{Laurent power series and first order thickening}
	Suppose $R\to R_0$ is a square-zero extension of $k$-algebras by $\mathfrak{a}\subset R$. Then for every field extension $K'/K$, the induced homomorphism $K'\otimes_K R\laurentseries{t}\to K'\otimes_K R_0\laurentseries{t}$ is a square-zero extension by the ideal $K'\otimes_K\mathfrak{a}\laurentseries{t}\subset K'\otimes_K R\laurentseries{t}$. \qed
\end{Lemma}

\begin{Lemma}\label{Lemma: Square-zero extension of torus for SU_3}
	Suppose $R\to R_0$ is a square-zero extension for $R_0,R\in\ArtLoc_k$. Then we have $\ker(LT_{SU_3}(R)\to LT_{SU_3}(R_0)) \subset LE_{SU_3}(R)$.
\end{Lemma}
\begin{proof}
	Put $\mathfrak{a}=\ker(R\to R_0)$. We have $\mathfrak{a}^2=0$. Then  \cref{Laurent power series and first order thickening} implies that the ideal $L\otimes_K\mathfrak{a}\laurentseries{t}\subset L\otimes_K R\laurentseries{t}$ is also square-zero. From the explicit description of $T_{SU_3}$ given in \cref{section:Kneser-Tits-for-SU3} we know that elements of $\ker(LT_{SU_3}(R)\to LT_{SU_3}(R_0))$ are of the form
	\[\diag(1+x,1-x+\sigma(x),1-\sigma(x))\]
	where $x\in L\otimes_K\mathfrak{a}\laurentseries{t}$. Then we have an explicit identity
	\[\begin{pmatrix}
		1 & x & 0 \\
		0 & 1 & \sigma(x)  \\
		0 & 0 & 1
	\end{pmatrix}
	\begin{pmatrix}
		1 & 0 & 0 \\
		1 & 1 & 0  \\
		z & 1 & 1
	\end{pmatrix}
	\begin{pmatrix}
		1 & -x & 0 \\
		0 & 1 & -\sigma(x)  \\
		0 & 0 & 1
	\end{pmatrix}
	=
	\begin{pmatrix}
		1+x & 0 & 0 \\
		* & 1-x+\sigma(x) & 0  \\
		* & * & 1-\sigma(x)
	\end{pmatrix}.\]
	One sees that the first and the third matrices on the left-hand side belong to $U_{SU_3}^+(R\laurentseries{t})$. The second matrix on the left-hand side belongs to $U_{SU_3}^-(R\laurentseries{t})$ provided that $z+\sigma(z)=1$. Such $z\in L$ exists by \cref{Lemma: Trace map is surjective}.
\end{proof}

\begin{Lemma}\label{lemma:identity-up-to-square-zero}
	Suppose $R\to R_0$ is a square-zero extension for $R_0,R\in\ArtLoc_k$. Then we have $\ker(LSU_3(R)\to LSU_3(R_0)) \subset LE_{SU_3}(R)$.
\end{Lemma}
\begin{proof}
	Pick $g\in \ker(LSU_3(R)\to LSU_3(R_0))$. Then it can be easily checked that the diagonal of $g$ is contained in $LT_{SU_3}(R)$. Using \cref{Lemma: Square-zero extension of torus for SU_3}, we may assume that the diagonal entries of $g$ are all equal to $1$. Put
	\[g=\begin{pmatrix}
		1 & a & b \\
		d & 1 & c  \\
		e & f & 1
	\end{pmatrix}, \quad
	u^+=\begin{pmatrix}
		1 & a & b \\
		0 & 1 & c  \\
		0 & 0 & 1
	\end{pmatrix}, \quad
	u^-=\begin{pmatrix}
		1 & 0 & 0 \\
		d & 1 & 0  \\
		e & f & 1
	\end{pmatrix}.\]
	\cref{Laurent power series and first order thickening} shows that  $a,b,c,d,e,f$ are contained in the square-zero ideal $L\otimes_K\mathfrak{a}\laurentseries{t}\subset L\otimes_K R\laurentseries{t}$, where $\mathfrak{a}=\ker(R\to R_0)$. From the explicit description of $SU_3$ given in \cref{section:Kneser-Tits-for-SU3}, we deduce that $u^+\in U_{SU_3}^+(R\laurentseries{t})$, $u^-\in U_{SU_3}^-(R\laurentseries{t})$ and $g=u^+u^-$.
\end{proof}

\begin{Lemma}\label{lemma:Artinian-Kneser-Tits-for-SU3}
	For any $R\in\ArtLoc_k$, we have $LT_{SU_3}(R)\subset LE_{SU_3}(R)$.
\end{Lemma}
\begin{proof}
	\cref{Field Kneser-Tits for SU_3} shows that the statement is true when $R$ is a field. Since Artinian local $k$-algebras can be obtained by successive square-zero extensions starting from a field, it suffices to show that the statement is true for every square-zero extension $R\to R_0$ provided that the statement is true for $R_0\in\ArtLoc_k$. 
	
	Pick an element $h\in LT_{SU_3}(R)$. Let $h'$ be the image of $h$ in $LT_{SU_3}(R_0)$. By assumption, we have $h'\in LE_{SU_3}(R_0)$. Hence we can find $u_1^+,\dots,u_m^+\in U_{SU_3}^+(R_0\laurentseries{t})$ and $u_1^-,\dots,u_m^-\in U_{SU_3}^-(R_0\laurentseries{t})$ such that $h'=u_1^+u_1^-\dots u_m^+u_m^-$. \cref{Laurent power series and first order thickening} implies that $R\laurentseries{t}\to R_0\laurentseries{t}$ is a square-zero extension. Since both $U_{SU_3}^+$ and $U_{SU_3}^-$ are smooth over $K$, the infinitesimal lifting property provides $v_1^+,\dots,v_m^+\in U_{SU_3}^+(R\laurentseries{t})$ and $v_1^-,\dots,v_m^-\in U_{SU_3}^-(R\laurentseries{t})$ which lift $u_1^+,\dots,u_m^+\in U_{SU_3}^+(R_0\laurentseries{t})$ and $u_1^-,\dots,u_m^-\in U_{SU_3}^-(R_0\laurentseries{t})$. Put
	\[h''=v_1^+v_1^-\dots v_m^+v_m^-\in LE_{SU_3}(R)\subset LSU_3(R).\]
	Then $h''$ and $h$ have the same image in $LSU_3(R_0)$. By \cref{lemma:identity-up-to-square-zero}, $h''h^{-1}$ is contained in the group $LE_{SU_3}(R)$, and hence so is $h$.
\end{proof}

\subsection{Type $A_2$}\label{section:type-A_2}
Assume the Dynkin scheme $\mathfrak{D}_i$ is of type $A_2$ for $i\in I$. Then we have $\Dyn(H_i)\cong\mathfrak{D}_i^{\min}=\mathfrak{D}_i$ by \cref{lemma:properties-of-root-groups}. We have shown that $H_i$ is quasi-split, semisimple and simply connected. The fact that $\mathfrak{D}_i$ is connected implies that $H_i$ is not split. Hence $H_i$ is isomorphic to the group $SU_3$ explicitly given in \cref{section:Kneser-Tits-for-SU3}, where the quadratic extension is $L=\mathcal{O}(\mathfrak{D}_i)$. As before, we can choose an isomorphism $H_i\simto SU_3$ which maps the Borel pair $T_i\subset B_i^+$ to $T_{SU_3}\subset B_{SU_3}^+$. Now \cref{lemma:Artinian-Kneser-Tits-for-SU3} implies the following statement.

\begin{Lemma}\label{lemma:Artinian-Kneser-Tits-for-type-A_2}
	If $\mathfrak{D}_i$ is of type $A_2$, then $LT_i(R)\subset LE_{H_i}(R)$ for every $R\in\ArtLoc_k$. \qed
\end{Lemma}

\subsection{Completion of the proof}
By definition, we have $U_i^+\subset U^+$ and $U_i^-\subset U^-$ for every $i\in I$. Hence $LE_{H_i}(R)\subset LE_G(R)$ for all $k$-algebras $R$. Now the isomorphism $\prod_{i\in I} T_i\simto T$ given by \cref{lemma:properties-of-root-groups} reduces \cref{proposition:Artinian Kneser-Tits for all G} to the following statement.

\begin{Proposition}
	We have $LT_i(R)\subset LE_{H_i}(R)$ for all $i\in I$ and $R\in\ArtLoc_k$.
\end{Proposition}
\begin{proof}
	The statement follows from \cref{lemma:type-of-connected-component-of-dynkin-scheme,lemma:Kneser-Tits-for-isotypic-of-type-A_1,lemma:Artinian-Kneser-Tits-for-type-A_2}.
\end{proof}

\bibliography{Reference}
\bibliographystyle{amsplain}

\end{document}